\documentclass[11pt]{0309}

\input nothing.d

\newif\ifnotesw\noteswtrue
\renewcommand{\comment}[1]{\ifnotesw $\blacktriangleright$\ {\sf #1}\ 
  $\blacktriangleleft$ \fi}
\noteswfalse	

\newcommand{\OK}{\mathrm{OK}}

\newcommand{\Cov}{\mathrm{Cov}}
\newcommand{\Var}{\mathrm{Var}}
\newcommand{\tower}{\mathrm{TOWER}}
\newcommand{\CG}{{\cal G}}
\newcommand{\CY}{{\cal S}} 
\newcommand{\Poisson}[1]{{\C P}_{#1}}
\newcommand{\ehr}[3]{\mbox{\sc Ehr}_{#1}(#2,#3)}



\newcommand{\ah}{\alpha}
\newcommand{\gam}{\gamma}

\newcommand{\ben}{\begin{enumerate}}
\newcommand{\een}{\end{enumerate}}
\newcommand{\bec}{\begin{center}}
\newcommand{\enc}{\end{center}}

\title{How Complex are Random Graphs in First Order Logic?}

\author{Jeong Han Kim\\
 Microsoft Research\\
 One Microsoft Way\\
 Redmond, WA 98052\\
 E-mail: \texttt{jehkim@microsoft.com}\\
 \mbox{}
 \and
Oleg Pikhurko\\
Department of Mathematical Sciences\\
Carnegie Mellon University\\
Pittsburgh, PA 15213-3890\\
Web: {\tt http://www.math.cmu.edu/\symbol{126}pikhurko}\\
\mbox{}
 \and
 Joel H.\ Spencer\\
Courant Institute\\ 
New York University\\
New York, NY 10012\\
E-mail: {\tt spencer@cs.nyu.edu}\\
\mbox{}
 \and
 Oleg Verbitsky\\
Department of Mechanics and Mathematics\\
Kyiv University\\ 
Ukraine\\
E-mail: {\tt oleg@ov.litech.net}\\
\mbox{}}

\date{November 17, 2003}
\maketitle

\begin{abstract}
 It is not hard to write a first order formula which is true for a
given graph $G$ but is false for any graph not isomorphic to $G$. The
smallest number $D(G)$ of nested quantifiers in a such formula can
serve as a measure for the ``first order complexity'' of $G$. 

Here, this parameter is studied for random graphs. We determine it
asymptotically when the edge probability $p$ is constant; in fact,
$D(G)$ is of order $\log n$ then. For very sparse graphs its magnitude
is $\Theta(n)$. On the other hand, for certain (carefully chosen)
values of $p$ the parameter $D(G)$ can drop down to the very slow
growing function $\log^* n$, the inverse of the $\tower$-function.
The general picture, however, is still a mystery.\end{abstract}

\section{Introduction}\label{intro}

In this paper we shall deal with sentences about graphs expressible in first order logic. Namely, the vocabulary consists of the following symbols:
 \begin{itemize}
 \item \emph{variables} ($x$, $y$, $y_1$, etc);
 \item the \emph{relations} $=$ (equality) and $\sim$ (the graph
adjacency);
  \item the \emph{quantifiers} $\forall$ (universality) and $\exists$
(existence);
 \item the usual Boolean \emph{connectives} ($\vee$,
$\wedge$, $\neg$, $\Leftrightarrow$, and $\Rightarrow$).
 \end{itemize}
 These can be combined into first order \emph{formulas} accordingly
to the standard rules.
 A \emph{sentence} is a formula without free variables. On the
intuitive level it is perfectly clear what we mean when we say that a
sentence \emph{is true} on a graph $G$. This is denoted by $G\models
A$; we write $G\not\models A$ for its negation ($A$ \emph{is false} on
$G$). We do not formalize these notions; a more detailed discussion
can be found in e.g.~\cite[Section~1]{spencer:slrg}. 

Please note that the variables represent vertices so the quantifiers
apply to vertices only, i.e.\ we cannot express sentences like
\textsl{``There is a set $X$ having a given property''}. (In fact, the
language lacks any symbols to represent sets or functions.) We do not
allow infinite sentences. As we do not go beyond first
order logic, the standalone term ``sentence'' means a ``first
order sentence.''

From a logician's point of view, the first order properties of graphs
form a natural class of properties to study. For example, the
so-called \emph{zero-one laws} for random graphs have been extensively
studied: see
e.g.~\cite{glebskii+:69,fagin:76,kolaitis+proemel+rothschild:85,kolaitis+proemel+rothschild:87,compton:87,shelah+spencer:88,luczak+spencer:91,shelah:96,spencer:slrg}.

Of course, if $G\models A$ and $H\cong G$ (i.e.\ $H$ is isomorphic to
$G$), then $H\models A$.  On the other hand, for any graph $G$ it is
possible to find a first order sentence $A$ which \emph{defines} $G$, that
is, $G\models A$ while $H\not\models A$ for any $H\not\cong
G$. Indeed, let $V(G)=\{v_1,\dots,v_n\}$. The required sentence $A$
could read:
 \begin{quote} \textsl{``There are vertices $x_1,\dots,x_n$, all
distinct, such that any vertex $x_{n+1}$ is equal to one of these and
$x_i\sim x_j$ iff $\{v_i,v_j\}\in E(G)$, $1\le i<j\le n$.''}
 \end{quote}

However, this sentence looks rather wasteful: we have $n+1$ variables,
the $\sim$-relation was used $n\choose 2$ times, etc. Of a number of
possible parameters measuring how complex $A$ is, we choose here
$D(A)$, the \emph{quantifier depth} (or simply \emph{depth}) which is
the size of a longest sequence of embedded quantifiers. (In the above
example, $D(A)=n+1$.) This is a natural characteristics which appears,
for example, in the analysis of algorithms for checking whether
$G\models A$. Also, the depth function can be studied by using the
so-called Ehrenfeucht game~\cite{ehrenfeucht:61} (see
Section~\ref{game} here). Following Pikhurko, Veith and
Verbitsky~\cite{pikhurko+veith+verbitsky:03} (see
also~\cite{pikhurko+verbitsky:03}) we let $D(G)$ be the smallest depth
of a sentence defining $G$. It is a measure of how difficult it is to
describe the graph $G$ in first order logic.

A word of warning: the function $D(G)$ does not correlate very well
with our everyday intuition of how complex the graph $G$ is. Such are
the limitations of the first order language that, for example,
$D(K_n)=D(\O K_n)=n+1$ is the largest among all order-$n$ graphs but
what can be simpler than the complete or empty graph?!

Perhaps, this is just an unlucky example? This approach seems only
partially helpful: it is shown in~\cite{pikhurko+veith+verbitsky:03}
that, certain exceptions aside, $D(G)$ can be bounded by
$\frac{v(G)+5}2$; however, the situation seems to get messy when we
try to describe all graphs with $D(G)\ge (\frac12-\e) v(G)$. Besides a
large homogeneous set, there are other obstacles which may push the
depth up: Cai, F\"urer, and Immerman~\cite{cai+furer+immerman:92}
constructed graphs $G$ to define which we need $\Omega(v(G))$ nested
quantifiers even if we add \emph{counting} to  first order
logic. This is a rather drastic addition: for example, we need only
two nested quantifiers to define $K_n$ with counting, namely,
 $$
 \mbox{\sl ``There are precisely $n$ vertices and every two of them are
   connected.''} 
 $$

The opposite approach was taken by Pikhurko, Spencer and
Verbitsky~\cite{pikhurko+spencer+verbitsky:03}: what is $g(n)$, the
minimum $D(G)$ over all graphs $G$ of order $n$? It turned out that
$g(n)$ can be arbitrarily small in the following sense: for any
recursive function $f:\I N\to\I N$ there is $n$ such that $f(g(n)) <
n$. If we try to ``smoothen'' $g(n)$ by defining $\gamma(n)=\max_{i\le
n} g(i)$, then $\gamma(n)=\Theta(\log^* n)$. Here, the \emph{log-star}
$\log^* n$ is the inverse to the $\tower$-function, that is, the
number of times we have to take the binary logarithm before we get
below one:
 $$
 \log^* n=\min\{i\in \I N\mid \log_2^{(i)} n<1\}.
 $$

Such a behavior is surprising and intriguing. Having studied the two
extreme cases, we concentrate now on what happens in a typical
graph. More generally, we consider the standard Erd\H os-R\'enyi model
$G\in \CG(n,p)$, where $p$ denotes the edge probability.  Of course,
we are interested in events occurring \emph{whp} (with high
probability, that is, with probability $1-o(1)$ as
$n\to\infty$). While a zero-one law studies the probability that a
fixed sentence holds, we take a random graph and ask what the
`simplest' sentence defining it is.

As $D(G)=D(\O G)$, we can assume without loss of generality that $p\le
\frac12$.

In Section~\ref{constant} we study the case when $0<p\le \frac12$ is a
constant and show that whp 
 \beq{intro:constant}
 D(G)=\log_{1/p} n + O(\ln \ln n).
 \eeq

The case $p=\frac12$ is always of particular interest:
$G\in\CG(n,\frac12)$ is uniformly distributed among all graphs of
order $n$. In Section~\ref{1/2}, we have found a different line of
argument (as far as the upper bound is concerned), which allowed us to
pinpoint $D(G)$ down to at most $5$ distinct values for infinitely
many $n$. Unfortunately, this approach does not seem to work
for $p\not=\frac12$.

In Section~\ref{sparse} we show that for $p<\frac{1.19}n$, $D(G)$ is
determined by the number of isolated vertices and therefore is of
order $\Theta(n)$. We believe (cf.\ Conjecture~\rcj{giant}) that the
giant component, which appears around $p=\frac1n$, has negligible
effect on $D(G)$ as long as $p=O(n^{-1})$ but we were not able to
prove this.

Rather surprisingly, for some carefully selected $p=p(n)$ the function
$D(G)$ can be as small as $O(\log^* n)$. The reason is that the
integer arithmetics can be modeled over the obtained random graphs
while integers can be defined by first order sentences of very
small depth. We do not present an exhaustive general theorem but give
an example demonstrating this phenomenon when $p=n^{-1/4}$. On the
other hand, the upper bound $O(\log^*n)$ is sharp, up to a
multiplicative constant, cf.\ Theorem~\rth{logstar}.

The first order complexity of $G\in\CG(n,p)$ for the general $p$
remains a mystery. Open problems and conjectures are scattered
throughout the text. See also Section~\ref{open} for some concluding
remarks.

\section{The Ehrenfeucht Game}\label{game}

For non-isomorphic graphs $G$ and $G'$ let $D(G,G')$ be the smallest
quantifier depth of a first order sentence $A$ \emph{distinguishing}
$G$ from $G'$ (that is, $G\models A$ while $G'\not\models A$). As the
negation sign does not affect the depth, we have $D(G,G')=D(G',G)$.

\blm{D} For any graph $G$ we have
 \beq{D}
 D(G)=\max\{ D(G,G')\mid G'\not \cong G\}.
 \eeq
 \elm
 \bpf Clearly, $D(G,G')\le v(G)+1$, so the right-hand side of~\req{D}
is well-defined. Theorem~2.2.1 in~\cite{spencer:slrg} implies that all
graphs can be split into finitely many classes so that any first order
sentence of depth at most $v(G)+1$ does not distinguish graphs within
a class. For each class, except the one which contains $G$, pick a
representative $G'$ and let $A_{G'}$ be a minimum depth sentence
distinguishing $G$ from $G'$. The disjunction of these $A_{G'}$ proves
the `$\le$'-inequality in~\req{D}.

The converse inequality is trivial.\epf

In the remainder of this section, we describe the Ehrenfeucht
game which is a very useful combinatorial tool for studying
$D(G,G')$. It was introduced by
Ehrenfeucht~\cite{ehrenfeucht:61}. Earlier, Fra\"\i
ss\'e~\cite{fraisse:54} suggested an essentially equivalent way to
compute $D(G,G')$ in terms of partial isomorphisms between $G$ and
$G'$. A detailed discussion of the game can be found
in~\cite[Section~2]{spencer:slrg}.

Let $G$ and $G'$ be two graphs. By replacing $G'$ with an isomorphic
graph, we can assume that $V(G)\cap V(G')=\emptyset$. The
\emph{Ehrenfeucht game} $\ehr kG{G'}$ is played by two players, called
\emph{Spoiler} and \emph{Duplicator} and consists of $k$ rounds. For
brevity, let us refer to Spoiler as `him' and to Duplicator as
`her'. In the $i$-th round, $i=1,\dots,k$, Spoiler selects one of the
graphs $G$ and $G'$ and marks some its vertex by $i$; Duplicator must
put the same label $i$ on a vertex in the other graph. (A vertex may
receive more than one mark.) At the end of the game (i.e.\ after $k$
rounds) let $x_1,\dots,x_k$ be the vertices of $G$ marked $1,\dots,k$
respectively, regardless of who put the label there; let
$x_1',\dots,x_k'$ be the corresponding vertices in $G'$. Duplicator
wins if the correspondence $x_i\leftrightarrow x_i'$ is a partial
isomorphism, that is, we require that $\{x_i,x_j\}\in E(G)$ iff
$\{x_i',x_j'\}\in E(G')$ as well as that $x_i=x_j$ iff
$x_i'=x_j'$. Otherwise, Spoiler wins.

The crucial relation is that for any non-isomorphic $G$ and $G'$ the
smallest $r$ such that Spoiler has a winning strategy in $\ehr rG{G'}$
is equal to $D(G,G')$. In fact, an explicit winning strategy for
Spoiler gives us an explicit sentence distinguishing $G$ from $G'$.

If Spoiler can win the game, alternating between the graphs $G$ and
$G'$ at most $r$ times, then the corresponding sentence has the
\emph{alternation number} at most $r$, that is, any chain of nested
quantifiers has at most $r$ changes between $\exists$ and
$\forall$. (To make this well-defined, we assume that no quantifier is
within the range of a negation sign.) Let $D_r(G)$ be the smallest
depth of a sentence which defines $G$ and has the alternation number
at most $r$. Clearly, $D_r(G)=\max\{D_r(G,G')\mid G'\not\cong G\}$,
where $D_r(G,G')$ may be defined as the smallest $k$ such that Spoiler
can win $\ehr kG{G'}$ with at most $r$ alternations.

For small $r$, this is a considerable restriction, giving a
qualitative strengthening of the obtained results. Therefore, we make
the extra effort of computing the alternation number given by our
strategies if the obtained $r$ is really small.

Finally, let us make a few remarks on our terminology. When a player
marks a vertex, we may also say that the player \emph{selects} (or
\emph{claims}) the vertex. Duplicator \emph{loses after $i$ rounds} if
the correspondence between $(x_1,\dots,x_i)$ and $(x_1',\dots,x_i')$
is not a partial isomorphism. (Of course, there is no point in
continuing the game in this situation.)

\section{Constant Edge Probability}\label{constant}

As $D(G)=D(\O G)$, we can assume without loss of generality that $p\le
\frac12$. For brevity let us denote $q=1-p$. In this section we prove
the following result.

\bth{constant} Let $p$ be a constant, $0<p\le \frac12$. Let
$G\in\CG(n,p)$. Then whp
 \beq{constant}
 - O(1) \le D(G) - \log_{1/p} n + 2\log_{1/p} \ln
n \le (2+o(1))\, \frac{\ln \ln n}{-p\ln
p -q \ln q}
 \eeq
 \eth

The lower bound follows by observing that if for any disjoint
$A,B\subset G$ with $|A|+|B|\le k$, there is a vertex $y$ connected to
everything in $A$ but to nothing in $B$ (this is called the
\emph{$k$-extension property} or the \emph{$k$-Alice Restaurant
property}), then $D(G)\ge k+2$.  The upper bound is obtained by some
kind of recursion, where for every $x\in G$ we write a sentence $A_x$
describing its neighborhood $\Gamma(x)$. Whp no two neighborhoods are
isomorphic so $A_x$ ``defines'' $x$ and the final sentence $A$
stipulates that the (unique) vertices satisfying $A_x$ and $A_y$ are
connected if and only if $x,y\in G$ are. As each recursive step
reduces the order by a factor of about $1/p$, the obtained sentence
has depth around $\log_{1/p}n$. There are some technicalities to
overcome. However, the combinatorial setting of the Ehrenfeucht games
makes the proof more transparent and accessible.

Unfortunately, we have hardly any control on the alternation number in
Theorem~\rth{constant}. The following result fills this gap by
providing the defining sentences of a very restrictive form: no
alternation at all. This, however, comes at the expense of increasing
the depth by a constant factor. 

\bth{D0} Let $p$ be a constant, $0<p<1$. Let
$G\in\CG(n,p)$. Then whp
 \beq{D0}
 D_0(G) \le (2 + o(1))\, \frac{\ln n}{-\ln(p^2+q^2)}.
 \eeq
 \eth

\brm If we are happy to bound $D_1$ only, then the constant
in~\req{D0} can be improved: in the proof (Section~\ref{D0}) we have
to use Lemma~\rlm{Y} instead of Lemma~\rlm{DetD0}.\erm

\subsection{The Lower Bound}

To prove the lower bound in~\req{constant} we use the following lemma.

\blm{alice} If $G$ has the $k$-extension property, then $D(G)\ge
k+2$.\elm
 \bpf Let $G'\not\cong G$ be another graph which has the $k$-extension
property. (For example, we can take a random graph of large order.)
Consider $\ehr {k+1}G{G'}$. Duplicator's strategy is
straightforward. If in the $i$-th round Spoiler selects a previously
marked vertex, Duplicator does the same in the other graph. Otherwise,
she matches the adjacencies between $x_i$ and $\{x_1,\dots,x_{i-1}\}$
to those between $x_i'$ and $\{x_1',\dots,x_{i-1}'\}$ by the
$k$-extension property.\epf
 
It is easy to show that whp $G\in\CG(n,p)$, for constant
$p\in(0,\frac12]$, has the $\floor{r}$-extension property with
 \beq{lower}
 r=\log_{1/p}n - 2 \log_{1/p} \ln n + \log_{1/p} \ln(1/p) -o(1),
 \eeq
 which gives us the required lower bound by Lemma~\rlm{alice}.
Indeed, for $k<r-\Theta(1)$, the expected number of `bad' $A,B\subset V(G)$
with $|A|+|B|=k$ can be bounded by
 $$
 {n\choose k} 2^k (1-p^k)^{n-k} = \me^{k\ln n -p^{k}n + o(\ln^2 n)}=o(1). 
 $$

\subsection{The Upper Bound}

Let $G=(V,E)$ be a graph. Let $\C V_i$ consist of all ordered
sequences of $i$ pairwise distinct vertices of $G$. For $\B
x=(x_1,\dots,x_i)\in\C V_i$ define
 $$
 V_{\B x}=\{y\in V\setminus\{x_1,\dots,x_i\} \mid \forall j\in[i]\
\{y,x_j\}\in E\},
 $$
 and $G_{\B x}=G[V_{\B x}]$. We abbreviate $G_{x_1,\dots,
x_i}=G_{(x_1,\dots,x_i)}$, etc. Let us agree that $\C
V_{-1}=\emptyset$, $\C V_0=\{()\}$ consists of the empty
sequence, and $G_{()}=G$. 
 
The following lemma specifies our global line of attack.

\blm{lupper} Suppose that a graph $G$, numbers $l\ge 0$ and $l_0\ge 3$
satisfy all of the following conditions.
 \begin{enumerate}
 \bit1 For any $\B x\in\C V_{l-1}\cup \C V_l$ we have $D(G_{\B x})\le
l_0$.
 \bit3 For any $i\le l-1$, $\B x\in\C V_i$, and distinct $y,z\in V_{\B x}$,
the following two conditions hold. Let $U=V_{\B x,y,z}$. 
 
 \begin{enumerate}
 \bit{3a} Any injection $f:U\to V_{\B x,y}$ which embeds $G_{\B x,y,z}$
as an induced subgraph into $G_{\B x,y}$ is the identity mapping. (In
particular, $G_{\B x,y,z}$ admits no non-trivial automorphism.)
 \bit{3b} There is a vertex $v\in V_{\B x,y}\setminus U$ such that for
any vertex $w\in V_{\B x,z}\setminus U$ we have $\Gamma(v)\cap
U\not=\Gamma(w)\cap U$, where $\Gamma$ denotes the neighborhood of
a vertex.
 \end{enumerate}
 \bit4 For any $i\le l-1$, $\B x\in \C V_i$, and distinct $y,z,w\in
V_{\B x}$, $G_{\B x,y,z}$ is not isomorphic to an induced subgraph of
$G_{\B x,w}$.
 \end{enumerate}
 Then $D(G)\le l+l_0$.
 \elm
 \bpf Let us observe first that Condition~\rit{3} (or
Condition~\rit{4}) implies that
 \beq{2}
 G_{x_1,\dots, x_i, y}\not\cong G_{x_1,\dots, x_i, z},\quad
\mbox{for any $i\le l-1$, $(x_1,\dots,x_i)\in\C V_i$, and distinct
$y,z\in V_{x_1,\dots,x_i}$.}
 \eeq

We prove the lemma by induction on $l$. If $l$ is $0$ or $1$,
then Condition~\rit1 alone implies the claim. So, let $l\ge 2$. Let
$G'=(V',E')$ be any graph which is not isomorphic to $G$.

\case1 Suppose that there is $x\in V$ such that $G_{x}\not\cong
G'_{y}$ for any $y\in V'$. Spoiler selects this $x$. Let $x'$ be
the Duplicator's reply. The graph $G_{x}$ satisfies all the assumptions of
Lemma~\rlm{lupper} with $l$ decreased by $1$. Spoiler will always play
inside one of $G_x$ or $G'_{x'}$. We can assume that Duplicator does
the same for otherwise the adjacencies to $x$ and $x'$ do not
correspond. As $G_x\not\cong G'_{x'}$, Spoiler can use the induction
to win the $(G_x,G'_{x'})$-game in at most $l+l_0-1$ moves, as
required.\ecpf\smallskip

The same argument works if there is $x\in V'$ such that $G'_x\not\cong
G_y$ for any $y\in V$. 

\case2 Suppose now that there are $x\in V$ and distinct $y',z'\in V'$
such that
 \beq{case21}
 G_x\cong G'_{y'}\cong G'_{z'}.
 \eeq
 Spoiler selects $y'\in V'$. Assume that Duplicator replies with
$y=x$, for otherwise $G_y\not\cong G'_{y'}$ by~\req{2} and
Spoiler proceeds as in Case~1. Now Spoiler selects $z'$; let $z\in V$
be the Duplicator's reply. We can assume that
 \beq{case22}
 G_{y,z}\cong G'_{y',z'},
 \eeq
 for otherwise Spoiler applies the inductive strategy to the
$(G_{y,z},G_{y',z'})$-game, where $l$ is reduced by $2$.

We show that Spoiler can win in at most $3$ extra moves now. Let
$U=V_{y,z}$ and $U'=V'_{y',z'}$. Spoiler selects the vertex $v\in
V_{y}\setminus U$ given by Condition~\rit{3b}. Let $v'\in
V'_{y'}\setminus U'$ be the Duplicator's reply. By Condition~\rit{3a}
and~\req{case21}--\req{case22} there is a bijection $g:V'_{y'}\to
V'_{z'}$ which is the identity on $U'$ and induces an isomorphism of
$G'_{y'}$ onto $G'_{z'}$. Spoiler selects $w'=g(v')$. Whatever the
reply $w\in G_z\setminus U$ of Duplicator is, $\Gamma(v)\cap
U\not=\Gamma(w)\cap U$. But in $G'$ we have $\Gamma(v')\cap
U'=\Gamma(w')\cap U'$. Spoiler can point this difference with one more
move into $U$. The total number of moves is $5\le l+l_0$, as
required.\ecpf

By~\req2, the only remaining case is the following.

\case3 Suppose that there is a bijection $g:V\to V'$ such that for any
$x\in V$ we have $G_x\cong G_{g(x)}$.

As $G\not\cong G'$, there are $y,z\in V$ such that $g$ does not
preserve the adjacency between $y$ and $z$. Spoiler selects $y$. We
can assume that Duplicator replies with $y'=g(y)$ for otherwise
Spoiler proceeds as in Case~1. Now, Spoiler selects $z$ to which
Duplicator is forced to reply with $z'\not=g(z)$. Assume that
$G_{y,z}\cong G_{y',z'}$ for otherwise Spoiler applies the inductive
strategy for $l-2$ to these graphs. But then $G_{y,z}$ is an induced
subgraph of $G_w$, where $w=g^{-1}(z')$, contradicting Condition~\rit4.
\epf

In order to finish the proof of Theorem~\rth{constant} we apply
Lemma~\rlm{lupper} to a random graph $G\in\CG(n,p)$.

Let $l=\log_{1/p}n -C \log_{1/p}\ln n-1$ and $l_0= C_0 \ln \ln n$ where $C$ and
$C_0$ are constants such that $C>2$ and $C_0(-p\ln p-q\ln q)> C$. Let
$m=np^{l+1}=\ln^C n$. Let $\e>0$ be a small constant. Let $n$ be
sufficiently large.

Let $\C V=\cup_{i=0}^{l+1} \C V_i$. Observe that $|\C V|\le
\me^{(1+\e) \log_{1/p}n\, \ln n}=\me^{O(\ln^2n)}$.

\blm{sizes} Whp for any $i\le l+1$ and $\B x\in\C V_i$ we have
 \beq{sizes}
 \left|\, |V_{\B x}| - p^in\,\right| \le \e p^i n.
 \eeq
 \elm
 \bpf Fix some $\B x\in \C V_i$. The size of $V_{\B x}$ has the
binomial distribution with parameters $(n-i,p^i)$. By Chernoff's
bound (\cite[Appendix~A]{alon+spencer:pm}), the probability $p'$ that
this $\B x$ violates~\req{sizes} is
 \beq{pprime}
 p'\le 2\me^{-\e^2np^i/3}\le 2\me^{-\e^2 m/3}= o(|\C V|^{-1}).
 \eeq
 Thus the expected number of `bad' $\B x$'s is $o(1)$, giving the
required.\epf

\blm{3} Whp Condition~\rit{3} holds.\elm
 \bpf Fix $i\le l-1$, $\B x\in\C V_i$, and $y,z\in V_{\B x}$. Let
$U=V_{\B x,y,z}$, $W=V_{\B x,y}$, $u=|U|$, and $w=|W|$.

First we deal with Condition~\rit{3a}. Take $j\in[1,u]$. Let $g$ be
any injection from $U$ into $W$ such that $|U_g|=j$, where $U_g=\{v\in
U\mid g(v)\not=v\}$ consists of the elements moved by $g$. Let the
same symbol $g$ denote also the induced action on edges. Let $E_g$
consist of those $e\in {U\choose 2}$ such that $g(e)\not=e$. It is not
hard to see that
 $$
 |E_g|\ge {j\choose 2}+j(u-j) - \frac j2 = j(u-j/2-1).
 $$

We can find a set $D\subset E_g$ of size at least $|E_g|/3$ such that
$D\cap g(D)=\emptyset$. We do so greedily: choose any $e\in E_g$, move
$e$ to $D$, and remove $g(e),g^{-1}(e)$ from $E_g$ (if they belong
there). The probability that, for all $e\in D$, the $2$-sets $e$ and
$g(e)$ are simultaneously edges or non-edges is
$(p^2+q^2)^{|D|}$ because these events are independent. This gives
an upper bound on the probability that $g$ induces an isomorphism.

Given $j$, there are at most ${u\choose j} w^j$ choices of $g$. The
sequence $(\B x,y)$ (or $(\B x,y,z)$) violates~\req{sizes} with
probability at most $p'$, where $p'$ is as in~\req{pprime}. Thus we
can bound the probability that $\B x,y,z$ violate Condition~\rit{3a}
by
 $$
 2p'+\sum_{j=1}^u {u\choose j} w^j (p^2+q^2)^{j(u-j/2-1)/3} < 2p'+
 (p^2+q^2)^{(\frac13-\e)m}.
 $$

Hence, the expected number of bad witnesses $\B x,y,z$ is at most
 $$
 |\C V|\, \left(2p'+ (p^2+q^2)^{(\frac13-\e)m}\right) =o(1),
 $$
 giving the required by Markov's inequality.

To estimate the probability that~\rit{3b} fails, fix some $v\in
W\setminus U$. The probability that some vertex of $V_{\B
x,z}\setminus U$ has the same neighborhood in $U$ as $v$ is at most
$(p^2+q^2)^u$. We have $u> (1-\e)m$ with probability at least
$1-p'$. Hence, $v$ does not satisfy Condition~\rit{3b} with
probability at most
 $$
 |V_{\B x,z}|\,\left(p'+(p^2+q^2)^{(1-\e)m}\right)=o(|\C V|),
 $$
 finishing the proof.\epf

Condition~\rit{4} is verified similarly to the argument of
Lemma~\rlm{3}. (The proof is, in a way, even easier because $|V_{\B
x,y,z}\setminus V_w|=\Omega(m)$ whp.) All that remains to check is
Condition~\rit{1}. To deal with it, we need another strategic
lemma. For a subset $X$ of vertices of $G=(V,E)$, define the
equivalence relation $\equiv_X$ on $V$, called the
\emph{$X$-similarity}, by $x\equiv_X y$ iff $x=y$ or $x,y\in
V\setminus X$ satisfy $\Gamma(x)\cap X = \Gamma(y)\cap X$. This is an
equivalence relation. Let
 $$
 \CY(X)=\big\{x\in V\mid \forall y\in V\
 (y\not=x\ \imply\ 
 y\not\equiv_X x)\big\} \supset X.
 $$
 The vertices in $\CY(X)$ are \emph{sifted out} by $X$ (that is, are
uniquely determined by their adjacencies to $X$). We call $X$ a
\emph{sieve} if $\CY(X)=V$.

\blm{Y} Let $X\subset V$. Define $Y=\CY(X)$. If $\CY(Y)=V$,
then $D_1(G)\le |X|+3$.\elm
 \bpf
 Let $G'\not\cong G$. First, Spoiler selects all of $X$. Let
$X'\subset V'$ be the Duplicator's reply. Assume that Duplicator has
not lost yet. For the notational simplicity let us identify $X$ and
$X'$ so that $V\cap V'=X=X'$ and our both graphs coincide on $X$.  Let
$Y'= \CY_{G'}(X)$.

It is not hard to see that Spoiler wins in at most two extra moves
unless the following holds. For any $y\in Y\setminus X$ there is a
$y'\in Y'\setminus X$ (and vice versa) such that $\Gamma(y)\cap
X=\Gamma(y')\cap X$. Moreover, this bijective correspondence between
$Y$ and $Y'$ induces an isomorphism between $G[Y]$ and
$G'[Y']$. 

Clearly, if Duplicator does not respect this correspondence, she loses
immediately. Therefore, we may identify $Y$ with $Y'$. Let
$Z=V\setminus Y$ and $Z'= V'\setminus Y$. Let $z\in Z$ and define
 $$
 W'_{z}=\big\{z'\in Z'\mid \Gamma(z')\cap Y = \Gamma(z)\cap
Y\big\}
 $$
 
If $W'_z=\emptyset$, Spoiler wins in at most two moves. First, he
selects $z$. Let Duplicator reply with $z'\in Z'$. As the
neighborhoods of $z,z'$ in $Y$ differ, Spoiler can highlight this by
picking a vertex of $Y$. If $|W'_z|\ge 2$, then Spoiler selects some
two vertices of $W'_z$ and wins with at most one more move, as
required.

Hence, we can assume that for any $z$ we have $W'_z=\{f(z)\}$ for some
$f(z)\in Z'$. It is easy to see that $f:Z\to Z'$ is in fact a
bijection (otherwise Spoiler wins in two moves). As $G\not\cong G'$,
the mapping $f$ does not preserve the adjacency relation between some
$y,z\in Z$. Now, Spoiler selects both $y$ and $z$. Duplicator cannot
respond with $f(y)$ and $f(z)$; by the definition of $f$ Spoiler can
win in one extra move.\epf

By Lemma~\rlm{Y}, to complete the proof of Theorem~\rth{constant} it
suffices to verify that whp for any $\B x\in\C V_{l-1}\cup \C V_l$
there is an $(l_0-3)$-set $X\subset V_{\B x}$ such that, with respect
to $H=G_{\B x}$, $\CY(Y)=U$, where $Y=\CY(X)$ and $U=V_{\B
x}$. Let $k=l_0-3$ and $u=|U|$. Fix any $X\in{U\choose k}$. With
probability at least $1-p'$ we have $up^2\le (1+\e)m$. Conditioned on
this, $G_{\B x}$ is still constructed by choosing its edges
independently. The probability that a vertex $y\in U\setminus X$
belongs to $Y$ is
 \beq{ti}
 \sum_{i=0}^k {k\choose i} p^iq^{k-i}\left(1-p^iq^{k-i}\right)^{u-k-1}
 \eeq
 We want to bound this probability from below. Let, for example, $i_0=
pk - k^{1/2}\ln k$. For $i_0\le i\le k$ we have $(1-p^iq^{k-i})^u\ge
1-\e$ by the definition of $C_0$. Chernoff's bound implies that
$\sum_{i=i_0}^{k} {k\choose i} p^iq^{k-i} > 1-\e$ as this sum
corresponds to the Binomial distribution with parameters
$(k,p)$. Hence, the expression~\req{ti} is at least $(1-\e)^2 > 1-2\e$
and the expectation
 $$
 E[\,|Y|\,] > (1-2\e)u.
 $$

We construct the martingale $Y_0,\dots,Y_{u-k}$, where we expose the
vertices of $U\setminus X$ one by one and $Y_i$ is the expectation of
$|Y|$ after $i$ vertices have been exposed. Changing edges incident to
a vertex, we cannot decrease or increase $|Y|$ more than by two. By
Azuma's inequality (\cite[Theorem~7.2.1]{alon+spencer:pm}), the
probability that $|Y|$ drops, say, below $(1-3\e)u$ is at most
$\me^{-\Omega(m)}=o(|\C V|)$.
 \comment{Joel, could you please check whether this application of
Azuma's inequality is correct?}%
 Whp each $Y$ has at least $(1-3\e)u$ elements. The following simple
lemma completes our quest.

\blm{YY} Let $\e=\e(p)>0$ be sufficiently small. Whp for any $\B x\in
\C V_{l-1}\cup \C V_l$, every set $Y\subset V_{\B x}$ of size at least
$(1-3\e)u$, $u=|V_{\B x}|$, is a sieve in $G_{\B x}$.\elm
 \bpf Let $\B x$ satisfy~\req{sizes}. The expected number of bad
triples $(Y,y,z)$ (that is, the distinct vertices $y,z\in V_{\B
x}\setminus Y$ have the same neighborhood in a set $Y$ of size at
least $(1-3\e)u$) is
 $$
 \sum_{i\le 3\e u} {u\choose i} u^2 (p^2+q^2)^{u-i} =o(|\C V|).
 $$

The claim follows from~\req{pprime}.\epf

\subsection{Games with no Alternations}\label{D0}

Following our standard scheme, we first specify a graph property which
ensures the desired bound on $D(G)$ and then show that a random graph
satisfies this property whp.

\blm{DetD0} Assume that in a graph $G=(V,E)$ we can find $X\subset V$
such that
 \begin{enumerate} 
 \bit{d1}
  $X$ is a sieve;
 \bit{d2} $G[X]$ has no nontrivial automorphism;
 \bit{d3} $G$ has no other induced subgraph isomorphic to $G[X]$.
\end{enumerate} 
 Then $D_0(G)\le |X|+2$.
 \elm
 \bpf
 Let $G'$ be an arbitrary graph non-isomorphic to $G$.  For some $G'$
Spoiler plays all the time in $G$, for others he plays all the time in
$G'$.

We first describe the strategy when Spoiler plays in $G$. Spoiler
selects all vertices in $X$.  Suppose that Duplicator managed to
establish $\phi\,:\,X\to X'$, a partial isomorphism from $G$ to $G'$,
where $X'\subset V'$ is the set of Duplicator's responses.  Denote
$Z=V(G)\setminus X$ and $Z'=V(G')\setminus X'$.  We call two vertices,
$v\in Z$ and $v'\in Z'$ \emph{$\phi$-similar\/} if the extension of
$\phi$ which takes $v$ to $v'$ is a partial isomorphism from $G$ to
$G'$.  Four cases are possible:

\case1
 The $\phi$-similarity is a one-to-one correspondence between $Z$ and $Z'$.

\case2
 There is $v\in Z$ without a $\phi$-similar counterpart in $Z'$.

\case3
 There is $v'\in Z'$ without a $\phi$-similar counterpart in $Z$.

\case4
 There are $v'_1,v'_2\in Z'$ with the same $\phi$-similar counterpart
in $Z$.

In Case~1 there are $v_1,v_2\in Z$ with adjacency different from the
adjacency between their $\phi$-similar counterparts in $Z'$.  Spoiler
selects $v_1$ and $v_2$ and wins. In Case 2 Spoiler wins by selecting
the vertex $v$. In Cases 3 and 4 Spoiler fails in this way but plays
differently from the very beginning.

Namely, if there exist $X'$ and a partial isomorphism $\phi:X\to X'$
such that Cases 3 or 4 occur, Spoiler begins with selecting all
vertices in $X'$. Duplicator is forced to reply in accordance with
$\phi$ due to the conditions assumed for $X$.  Then Spoiler selects
the vertex $v'$ in Case 3 or $v'_1$ and $v'_2$ in Case 4 and wins.\epf

\blm{RandomD0} Let $\e>0$ and $0<p<1$ be fixed. Let
$G\in\CG(n,p)$ and let $X\subset V$ be any set of size $t\ge
(2+\e)\log_{1/r} n$, where $r=p^2+q^2$. Then whp
Conditions~\rit{d1}--\rit{d3} of Lemma~\rlm{RandomD0} hold.
 \elm
 \bpf
 The expected number of vertices with the same neighborhood in $X$ is
at most $n^2r^t=o(1)$, implying~\rit{d1}. Conditions~\rit{d2}
and~\rit{d3} follow from the following claim. 

\claim1 Whp no injective $g:X\to V$, with the exception of the
identity mapping, preserves the adjacency relation.
 \bcpf
 Fix $g$. Let $k=|K|$ and $l=|L|$, where
 \begin{eqnarray*}
 K&=&\{x\in X\mid g(x)\not\in X\},\\
 L&=&\{x\in X\mid g(x)\in X\setminus\{x\}\},
 \end{eqnarray*}
 As in the proof of Lemma~\rlm{3} we can find a set $D\subset
{X\setminus K\choose 2}$ of size at least $l\,\frac{t-k-l/2-1}3$ such
that $D\cap g(D)=\emptyset$. The latter property still holds if we
enlarge $D$ by the set of all elements of ${X\choose 2}$ incident to
$K$. Hence, the total probability of failure is at most
 $$
 \sum_{0<k+l\le t} {t\choose k} {t-k\choose l} n^k t^l
r^{l\,\frac{t-k-l/2-1}3 + k(t-1) - {k\choose 2}}=
 o(1),
 $$
 \comment{ Let $s_{k,l}$ denote the summand. First observe that both
$s_{1,0}=nr^{t}$ and $s_{0,1}=t^2r^{t/3}$ are $o(t^{-2})$. (I ignore
less significant terms.) We claim that for all other $(k,l)$, we have
$s_{k,l}=o(t^{-2})$, which would give the required.

Let $c>0$ be a small constant.

Suppose first that $k>(1-c)t$. Then $tk-k^2/2 > (1/2-c)t^2$, and
 $$
 s_{k,l}< n^{t} r^{(1/2-c)t^2} < n^t n^{-(2+\e)(1/2-c)t} = o(t^{-2}).
 $$
 (Here is the crucial place, where the restriction on $t$ comes into
play.)

So we can assume that $k<(1-c)t$.

If $l>ct$, then replacing $l$ by, say, $2$, we increase $s_{k,l}$: the
old/new ratio is $< t^t r^{ct (1-(1-c)-c/2)t/3}=o(1)$.

So assume that $l<ct$. Fixing $t$, let us decrease $l$ by one by one until
we hit $0$ (or $1$ if $k=0$). Each time the ratio is
 $$
 s_{k,l+1}/s_{k,l} < t^2 r^{\Omega(t)} =o(1),
 $$
 so $s_{k,l}$ increases. (Here we used the fact that $t-k-l/2=\Omega(t)$
as $k<(1-c)t$ and $k+l\le t$.)

If $(k,l)$ is $(1,0)$ or $(0,1)$, then we are done. So we can assume that
$l=0$ and $k\ge2$. Roughly, $s_{k,0}=n^kr^{kt-k^2/2}$. Let $k=xt$, then
this term is $n^{xt}n^{(2+\e)(x-x^2/2)t}$. But $x<2(x-x^2/2)$ for $0<x<1$,
so $s_{k,l}$ is still small unless $k=o(t)$. But then $s_{k,l}= n^kr^{kt}$
while by the choice of $t$ we have $nr^t<1/n$ and we are done.}%
 completing the proof of the claim and the lemma.\epf

\section{Edge Probability $1/2$}\label{1/2}

Here is the main result of this section. 

\bth{1/2} Let $G\in\CG(n,\frac12)$. For infinitely many values of
$n$ we have whp
 \beq{1/2}
 D_2(G)\le \log_2 n - 2\log_2 \ln n + \log_2 \ln 2+6+o(1).
 \eeq
 \eth

\brm The lower bound given by Lemma~\rlm{alice} and the case
$p=\frac12$ of~\req{lower} is by at most $5+o(1)$ smaller than the
upper bound in~\req{1/2}. This implies that $D(G)$ and $D_2(G)$ are
concentrated on at most $6$ different valued for such $n$. In
Section~\ref{together} we will show that whp we have only $5$
possible values.\erm

Before we start proving Theorem~\rth{1/2}, let us observe that for
$p=\frac12$ and an \emph{arbitrary} $n$ the upper
bound~\req{constant} can be improved by using Lemma~\rlm{Y} to
 $$
 D_1(G)\le \log_2 n - \log_2 \ln n +O(\ln\ln\ln n).
 $$
 (Details are left to the interested Reader.)

\subsection{Spoiler's Strategy}

Before we can specify the plan of our attack on Theorem~\rth{1/2}, we
have to give a few definitions.

Let $G=(V,E)$, $W\subset V$, and $u\in\I N$. Building upon the notions
defined before Lemma~\rlm{Y}, let
 $$
 \CY_u(W)= \bigcup_{U\subset V\setminus W\atop |U|=u} \Big(\CY(U\cup
 W)\setminus (U\cup W)\Big).
 $$
 In other words, a vertex $y\not\in W$ belongs to $\CY_u(W)$ iff there
is a $u$-set $U\subset V\setminus (W\cup\{y\})$ such that $U\cup W$
sifts out $y$. Note that $\CY(W)=\CY_0(W)\cup W$.

\blm{1/2} Let $Y=\CY_u(W)$. Suppose that $Y\cup W$ is a sieve in $G$
and that no two vertices of $Y$ have the same neighborhood in $W$. Then
$D_2(G)\le u+w+4$, where $w=|W|$.\elm
\bpf Let $G'=(V',E')$ be a graph non-isomorphic to $G$.
We describe a strategy allowing Spoiler to win the game $\ehr{u+w+4}G{G'}$.

Spoiler first claims $W$. Let Duplicator reply with $W'\subset V'$.
Assume that she does not lose in this phase, establishing a partial
isomorphism $f:W\to W'$. Recall that we call two vertices
$v\in V$ and $v'\in V'$ \emph{$f$-similar} if the extension
of $f$ taking $v$ to $v'$ is still a partial isomorphism from
$G$ to $G'$. Let $Y'=\CY_u(W')$.

\claim1
As soon as Spoiler moves inside $Y\cup Y'$ but Duplicator replies
outside $Y\cup Y'$, Spoiler is able to win in the next $u+1$ moves
with 1 alternation between the graphs.

\bcpf
Assume for example that, while Spoiler selects $y\in Y$, Duplicator
replies with $y'\notin Y'$. Spoiler selects some $u$-set $U$ with $y$
sifted out by $U\cup W$. Let the reply to it be $U'$. By the assumption on
$y'$, there is another vertex $z'$ with the same adjacencies to $U'\cup
W'$. Spoiler selects $z'$ and wins.\ecpf

\claim2
If $Y'$ contains two vertices with the same adjacencies to $W'$,
then Spoiler is able to win the game in $w+u+3$ moves with
at most 2 alternations.

\bcpf Assume that $y'$ and $z'$ are both in $Y'$ and have the same
adjacencies to $W'$. Spoiler selects these two vertices.  In order not
to lose immediately, Duplicator is forced to reply at least once
outside $Y$. Spoiler wins in the next $u+1$ moves according to
Claim~1.\ecpf

Assume therefore that all vertices in $Y'$ have pairwise distinct
neighborhoods in $W'$. This assumption and Claim 1 imply that
either the $f$-similarity determines a one-to-one correspondence
between $Y$ and $Y'$ or Spoiler is able to win the game in
$w+u+3$ moves with at most 2 alternations.
We will assume the first alternative.
Extend $f$ to a map from $W\cup Y$ onto $W'\cup Y'$ accordingly to
the $f$-similarity correspondence between $Y$ and~$Y'$.

\claim3
Suppose that Duplicator failed to respect the bijection $f$
after a Spoiler's move into $Y\cup Y'$. Then Spoiler can win in at most
$u+1$ extra moves, during which he alternates at most once.

\bcpf
Suppose that the previous move $x$ of Spoiler was in $G$, for
example. Clearly, the Duplicator's response $x'$ cannot belong to $Y'$
because $f(x)$ is the only vertex in $Y'$ with the required
$W'$-adjacencies. Spoiler applies the strategy of Claim~1.\ecpf

\claim4
If $f:W\cup Y\to W'\cup Y'$ is not a partial isomorphism from $G$ to $G'$,
then Spoiler is able to win the game in $w+u+3$ moves with 1 alternation.

\bcpf
Assume, for example, that $\{y_1,y_2\} \in {Y\choose 2}$ is an edge while
$\{f(y_1),f(y_2)\}$ is not. Spoiler picks $y_1$ and $y_2$.
Duplicator cannot reply with $f(y_1)$ and $f(y_2)$ so Spoiler
wins in at most $w + 2 + (u+1)= w+u+3$ moves by Claim~3.\ecpf

Assume therefore that $f:W\cup Y\to W'\cup Y'$ is a partial isomorphism.
Denote $R=V\setminus (W\cup Y)$ and $R'=V'\setminus (W'\cup Y')$.

\claim5
As soon as Spoiler moves inside $R\cup R'$ but Duplicator fails to reply
with an $f$-similar vertex in $R\cup R'$ (in the other graph),
Spoiler can win in at most
$u+2$ extra moves, during which he alternates at most once.

\bcpf
If Duplicator replies with a vertex $x\in Y\cup Y'$, in the next move
Spoiler marks the $f$-mate of $x$ and then applies the
strategy of Claim 3. If she replies in $R\cup R'$ but not with
an $f$-similar vertex, Spoiler highlights this in one more
move and again uses Claim~3.\ecpf

\claim6
If $W'\cup Y'$ is not a sieve in $G'$, then Spoiler
is able to win the game in $w+u+4$ moves with 2 alternations.

\bcpf
Spoiler picks two witnesses $z_1',z_2'\in R'$ with the same
adjacencies to $W'\cup Y'$. If at least one of the
corresponding replies $z_1,z_2$ is not in $R$, Spoiler
applies the strategy of Claim 5. Otherwise $z_1$ and $z_2$
belong to different $W\cup Y$-similarity classes
and there is a vertex $x\in W\cup Y$ adjacent to exactly one of
$z_1$ and $z_2$. If $x\in W$, Spoilers wins immediately.
If $x\in Y$, Spoiler picks $f(x)\in Y'$.
Duplicator cannot respond with $x$.
Now Spoiler wins the game in at most $u+1$ extra moves by Claim~3.\ecpf

In the rest of the proof we suppose that $W'\cup Y'$ is a sieve.
This assumption and Claim 5 imply that
either the $f$-similarity determines a one-to-one correspondence
between $R$ and $R'$ or Spoiler is able to win the game in
$w+u+3$ moves with 2 alternations. Let us assume the first alternative.
Extend $f$ to the whole of $V$ accordingly to the $f$-similarity
correspondence between $R$ and $R'$.
Thus $f$ is a bijection between $V$ and $V'$ now.
As $G$ and $G'$ are not isomorphic, $f$ does not preserve the
adjacency for some $\{y_1,y_2\}\in {V\choose 2}$. Spoiler selects
$y_1$ and $y_2$. If Duplicator replies with $f(y_1)$ and
$f(y_2)$, she loses immediately.
Otherwise, Spoiler applies the strategy of Claim 5
and wins, having made totally at most
$u+w+4$ moves and 1 alternation.
\epf

\subsection{The Probabilistic Part}

Let $k$ be given. For simplicity let us assume that $k$ is
even. Define
 $$
  f(n,k)={{n-\frac{k}{2}}\choose \frac{k}{2}}(n-k)(1-2^{-k})^{n-k-1}.
 $$
 
Basic asymptotics show that for $n=\Theta(k^2\,2^k)$ we have
$\frac{f(n+1,k)}{f(n,k)} \approx 1$ and thus we can find
$n=(\frac{\ln 2}2+o(1))\, k^2\,2^k$ such that
 $$
 f(n,k)=(10+o(1))\, \log_2 n.
 $$
 We fix this $n$. Routine calculations show that
 \beq{k}
  k \le \log_2 n -2 \log_2 \ln n + \log_2 \ln 2 + 1 +o(1).
 \eeq

Let $A$ be a fixed $\frac k2$-subset of $G\in\CG(n,\frac12)$. Let $\C
U$ consist of pairs $(U,y)$, where $U$ is a $k$-set containing $A$ and
$y\in V\setminus U$. For $(U,y)\in\C U$, let $I(U,y)$ denote the
indicator random variable for the event $y\in \CY(U)$.
We define
 \begin{eqnarray*}
 X &= & \sum_{(U,y)\in \C U} I(U,y),\\
 M &=& |\C U|\ =\ {n-k/2\choose k/2}(n-k),\\
 p &=& E[\,I(U,y)\,]\ =\ (1-2^{-k})^{n-k-1},\\
 \end{eqnarray*}
 We further set
 \begin{equation}\label{eq:mu}
 \mu = E[X] = Mp= f(n,k) = (10+o(1))\, \log_2 n. 
 \eeq

The idea behind these definitions is that we try to apply
Lemma~\rlm{1/2} for $W=A$ and $u=\frac k2$. Then $\mu$ is the expected
number of ways to construct a vertex of $\CY_u(W)$. Our proof works
only if $\mu$ is neither too big nor too small, that is, for some
special values of $n$ only. We do not know if $D(G)$ can pinned down
to $O(1)$ distinct values for an arbitrary $n$.

As $k\approx \log_2 n$ we have 
 $$
 M \approx \frac{n^{\frac{k}{2}+1}}{(k/2)!} = n^{\frac{k}{2}(1+o(1))}.
 $$
 As $\mu=n^{o(1)}$ we further have
 \begin{equation}\label{a}
 p\approx e^{-n2^{-k}} = n^{-\frac{k}{2}(1+o(1))}.
 \eeq

\blm{U1U2} For distinct $(U_1,y_1), (U_2,y_2)\in\C U$,
 \begin{equation}\label{f}
 E[\,I(U_1,y_1)\;I(U_2,y_2)\,] < n^{-\frac{3k}{4}(1+o(1))}.
 \eeq
 When $|U_1\cap U_2| < \frac{9k}{10}$
\begin{equation}\label{c}
 E[\,I(U_1,y_1)\;I(U_2,y_2)\,] < E[I(U_1,y_1)]\;E[I(U_2,y_2)]\;
(1+ O(n^{-\frac{k}{10}(1+o(1))})) \eeq
 \elm
 \bpf Condition on the adjacency patterns of $y_1,y_2$ to $U_1,U_2$
respectively.  
Let $z$ be any vertex not in $U_1\cup U_2\cup \{y_1,y_2\}$.
Suppose (the main case) $U_1\neq U_2$.  Then
 \begin{equation}\label{both} 
 \Pr[\,(z\equiv_{U_1}y_1) \wedge (z\equiv_{U_2}y_2)\,] \leq 2^{-k-1} \eeq
 as the adjacency pattern of $z$ to $U_1\cup U_2$ is then determined.
When $U_1=U_2$ the adjacency patterns of $y_1,y_2$ to $U_1$ must be
different as otherwise $I(U_1,y_1)=I(U_2,y_2)=0$.  Then it would
be impossible to have $z\equiv_{U_1}y_1$ and
$z\equiv_{U_2}y_2$ so (\ref{both}) still holds.
By inclusion-exclusion
 $$
 \Pr[\,(z\equiv_{U_1}y_1) \vee (z\equiv_{U_2}y_2)\,] \geq 
2\cdot 2^{-k} - 2^{-k-1} = 3\cdot 2^{-k-1}.
 $$
 If $I(U_1,y_1)\; I(U_2,y_2)=1$ then this fails for all such $z$.
But these events are mutually independent.  Thus, by (\ref{a})
 \begin{equation}\label{b}
 E[\,I(U_1,y_1)\;I(U_2,y_2)\,]\leq (1-3\cdot 2^{-k-1})^{n-2k-2}
= n^{-\frac{3k}{4}(1+o(1))}, \eeq
 giving the required.

Suppose further that $|U_1\cap U_2| < \frac{9k}{10}$.  Again let
$z$ be any vertex not in $U_1\cup U_2\cup \{y_1,y_2\}$ and
condition on the adjacency patterns of $y_1,y_2$ to $U_1,U_2$
respectively.  Now
 $$
 \Pr[\,(z\equiv_{U_1}y_1) \wedge (z\equiv_{U_2}y_2)\,] \leq
2^{-\frac{11k}{10}}
 $$
 as this event requires $z$ to have a given adjacency pattern to $U_1\cup U_2$.
Again by inclusion-exclusion
 $$
 \Pr[\,(z\equiv_{U_1}y_1) \vee (z\equiv_{U_2}y_2)\,] \geq 
 2\cdot 2^{-k} - 2^{-\frac{11k}{10}}.
 $$
 As with (\ref{b}) we deduce
 $$
 E[\,I(U_1,y_1)\;I(U_2,y_2)\,]\leq \Big(1-2\cdot 2^{-k}
+ 2^{-\frac{11k}{10}}\Big)^{n-2k-2}.
 $$
 Now we want to compare this to $E[I(U_1,y_1)]E[I(U_2,y_2)]$. We have
 $$
 (1-2\cdot 2^{-k}+2^{-\frac{11k}{10}})^{n-2k-2} =
((1- 2^{-k})^{n-k-1})^2(1+O(n\,2^{-\frac{11k}{10}})),
 $$
 yielding (\ref{c}).\epf

\blm{var}
 \begin{equation}\label{var} \Var[X] = O(E[X]). \eeq
 \elm
 \bpf As $X=\sum I(U,y)$, the sum of indicator random
variables, we employ the general bound
 \begin{equation}\label{e} \Var[X] \leq E[X] + \sum
\Cov[\,I(U_1,y_1),\,I(U_2,y_2)\,], \eeq
 the sum over distinct $(U_1,y_1),(U_2,y_2)\in \C U$.  The first term
is $\mu$.  Consider the sum of the covariances satisfying $|U_1\cap
U_2|> \frac{9k}{10}$.  There are $M$ choices for $(U_1,y_1)$.  For a
given $U_1$ there are $n^{\frac{k}{10}(1+o(1))}$ choices for
$(U_2,y_2)$ and
 \begin{equation}\label{w1} 
 \sum E[\,I(U_1,y_1)\;I(U_2,y_2)\,] \leq
Mn^{\frac{k}{10}(1+o(1))}n^{-\frac{3k}{4}(1+o(1))} = o(1). 
\eeq
 As the covariance of indicator random variables 
is at most the expectation of the product,
\begin{equation}\label{w7} \sum \Cov[I(U_1,y_1),I(U_2,y_2)] = o(1),  \eeq
 where in (\ref{w1})--(\ref{w7}) the sum is restricted to 
$|U_1\cap U_2|> \frac{9k}{10}$.  
When $|U_1\cap U_2|\leq \frac{9k}{10}$  we have from (\ref{c}) that
 $$
  \Cov[\,I(U_1,y_1),\,I(U_2,y_2)\,] =
O\left( E[I(U_1,y_1)]\;E[I(U_2,y_2)]\;
n^{-\frac{k}{10}(1+o(1))}\right).$$
  Hence (the sum over $|U_1\cap U_2|\leq \frac{9k}{10}$)  
 $$
 \sum \Cov[\,I(U_1,y_1),\,I(U_2,y_2)\,] =
O\Big(n^{-\frac{k}{10}(1+o(1))}\Big) 
\sum E[I(U_1,y_1)]\; E[I(U_2,y_2)]
 $$
 But $\sum E[I(U_1,y_1)]\;E[I(U_2,y_2)]$ over {\em all}
$(U_1,y_1),(U_2,y_2)\in \C U$ is precisely $\mu^2 = O(\ln^2n)$.  This
becomes absorbed in the $n^{-\frac{k}{10}}$ term and (the sum again
over $|U_1\cap U_2|\leq \frac{9k}{10}$)
 \begin{equation}\label{c3} \sum \Cov[\,I(U_1,y_1),\,I(U_2,y_2)\,] =
O(n^{-\frac{k}{10}(1+o(1))}) . \eeq
 In particular, all covariances in (\ref{e}) together add up to $o(1)$
and so we actually have the stronger result $\Var[X] \leq E[X] +
o(1)$.\epf

\blm{6} Whp every pair
$(U_1,y_1)\neq (U_2,y_2)$ from $\C U$ with $I(U_1,y_1)=I(U_2,y_2)=1$
satisfies
 \begin{enumerate}
 \bit{z1} $U_1\cap U_2 = A$;
 \bit{z2} $y_2\not\in U_1$;
 \bit{z4} $y_1\neq y_2$;
 \bit{z3} $y_1\not\equiv_A y_2$;
 \bit{z5} For $u_1,u_2\in U_1$ with $u_1\neq u_2$, $u_1\not\equiv_A u_2$;
 \bit{z6} For $u_1\in U_1, u_2\in U_2$, $u_1\not\equiv_A u_2$.
\end{enumerate}
 \elm
 \bpf From (\ref{w1}) the expected number of pairs
$(U_1,y_1)\neq (U_2,y_2)$ with $I(U_1,y_1)=I(U_2,y_2)=1$
and $|U_1\cap U_2|> \frac{9k}{10}$ is $o(1)$.  The total number
of pairs 
$(U_1,y_1)\neq (U_2,y_2)$ with $(U_1\cap U_2)\setminus A\neq\emptyset$ is less than 
$M^2\frac{k^2}{n}$.  For each with
$|U_1\cap U_2|\leq \frac{9k}{10}$ a weak form of (\ref{c}) gives that
$E[\,I(U_1,y_1)\;I(U_2,y_2)\,]\leq 2p^2$.  Hence the expected number of
such pairs with $I(U_1,y_1)=I(U_2,y_2)=1$ is bounded from above
by $M^2\frac{k^2}{n}(2p^2) = O((\ln^4n)/n)=o(1)$.  Hence the
probability that \rit{z1} fails is $o(1)$.

For \rit{z2} we first employ \rit{z1} and restrict attention
to $U_1\cap U_2 = A$. The number of $(U_1,y_1), (U_2,y_2)$ with
$y_2\in U_1$ is less than $M^2\frac{k}{n-k}$ and for each
$E[\,I(U_1,y_1)\;I(U_2,y_2)\,]\approx p^2$ so the expected number with
$I(U_1,y_1)\;I(U_2,y_2)=1$ is less than $\approx (Mp)^2\frac{k}{n}$ which
is $O((\ln^3n)/n)=o(1)$.

For \rit{z4} we first employ \rit{z1} and restrict attention
to $U_1\cap U_2 = A$. The number of such $(U_1,y_1),(U_2,y_2)$ with
$y_1=y_2$ is about $M^2n^{-1}$ and for each such 
$\Pr[I(U_1,y_1)=I(U_2,y_2)=1]\approx p^2$ so the expected number of violations
of \rit{z4} is around $(Mp)^2n^{-1} = \mu^2n^{-1}=o(1)$.

For \rit{z3} we first employ \rit{z1}--\rit{z4} and restrict attention
to pairs satisfying those conditions. The number of such pairs is
$\approx M^2$. For each such $\Pr[I(U_1,y_1)=I(U_2,y_2)=1]\approx p^2$
and, conditioned on this, $\Pr[y_1\equiv_A y_2] \approx 2^{-k/2}$.
Hence the expected number of pairs violating \rit{z3} is $\approx
(Mp)^22^{-k/2}=\mu^22^{-k/2}$ which is again $o(1)$.

For \rit{z5} there are $M$ choices of $U_1,y_1$ and $O(k^2)$ choices
for $u_1,u_2$.  For each
$\Pr[I(U_1,y_1)=1]=p$, $\Pr[u_1\equiv_A u_2] = 2^{-k/2}$ and these events
are independent so the expected number of violations of \rit{z5} is
$O(Mpk^22^{-k/2})$ which is $o(1)$.

For \rit{z6} we restrict attention to those cases satisfying
\rit{z1}--\rit{z3}.
There are less than
$M^2$ choices of $U_1,y_1,U_2,y_2$ and $O(k^2)$ choices
for $u_1,u_2$.  For each
$\Pr[I(U_1,y_1)=I(U_2,y_2)=1]\approx p^2$,
$\Pr[u_1\equiv_A u_2] = 2^{-k/2}$ and these events
are independent so the expected number of violations of \rit{z6} is
$O(M^2p^2k^22^{-k/2})$ which is $o(1)$.\epf

\blm{AlmostThere}
 Let $G\in \CG(n,\frac{1}{2})$ and $A$ be a fixed subset of $\frac{k}
{2}$ vertices, as above.  Let $Z$ denote the union of all sets $U-A$
where some $I(U,y)=1$ and let $S=\CY_{k/2}(A)$ denote the set of all
such vertices $y$.  Let $R=V\setminus(A\cup Z\cup S)$.  Then whp
 \begin{enumerate}
\bit{t1} All $y_1,y_2\in S$ have $y_1\not\equiv_A y_2$.
\bit{t2} All $u_1,u_2\in Z$ have $u_1\not\equiv_A u_2$.
\bit{t3} There are no distinct $z_1,z_2,z_3,z_4\in R$ with $z_1\equiv_S z_2$ and
$z_3\equiv_S z_4$.
\bit{t4} There are no distinct $z_1,z_2,z_3\in R$ with $z_1\equiv_S z_2 \equiv_S z_3$ 
\end{enumerate}
 \elm
 \comment{We believe that \rit{t3} and \rit{t4} could be replaced by the
stronger statement that there no not exist any distinct $z_1,z_2\in R$ with
$z_1\equiv_S z_2$ but our proof technique is not strong enough to show this.}%
 \bpf The first two statements are the conclusions \rit{z3}--\rit{z6} of
Lemma~\rlm{6}. We concentrate on showing \rit{t3} as \rit{t4} is
similar.  Set
 \newcommand{\LL}{l}
 \begin{equation}\label{L}
 \LL=\mu - \mu^{0.6}  \eeq
Let $Y$ denote the number of $\LL$-sets $\{(U_i,y_i)\mid 1\leq i\leq \LL\}$ 
(counting permutations of the $(U_i,y_i)$ as the same) and
$z_1,z_2,z_3,z_4$ satisfying
 \begin{itemize}
 \item $I(U_i,y_i)=1$ for $1\leq i\leq \LL$.
 \item The $U_i-A$ are disjoint, the $y_i$ are distinct, and no $y_i\in U_j$.
 \item $z_1,z_2,z_3,z_4\in R$ where $R$ denotes all vertices except the $U_i$ and the $y_i$
 \item $z_1\equiv_S z_2$ and $z_3\equiv_S z_4$ where we set $S=\{y_1,\ldots,y_\LL\}$.
\end{itemize}

We bound $E[Y]$.  There are less than $M^\LL/\LL!$ choices for the $(U_i,y_i)$ and
$n^4$ choices for the $z_1,z_2,z_3,z_4$.  Fix those choices. Set
$R^-=R\setminus\{z_1,z_2,z_3,z_4\}$ and
let $z\in R^-$.  For each $1\leq i\leq \LL$
 $$
 \Pr[z\equiv_{U_i}y_i] = 2^{-k}.
 $$
 For each $1\leq i,j\leq \LL$, as $|U_i\cap U_j|=\frac{k}{2}$,
 $$
 \Pr[\,(z\equiv_{U_i}y_i) \wedge (z\equiv_{U_j}y_j)\,] \leq 2^{-\frac{3k}{2}}.
 $$
 We apply the Bonferroni inequality, in the form that the probability of a disjunction
is at least the sum of the probabilities minus the sum of the pairwise probabilities:
 $$
 \Pr\left[\vee_{i=1}^\LL z\equiv_{U_i}y_i\right] \geq \LL\,2^{-k}-{\LL\choose 2}\,
2^{-\frac{3k}{2}}.
 $$
These events are independent over the $z\in R^-$ as they involve different adjacencies.
Let $\OK$ denote the event that no $z\equiv_{U_i}y_i$ for any $z\in R^-$ and
$1\leq i\leq \LL$. The independence gives:
 $$
 \Pr[\,\OK\,] \leq
\left(1-\LL\,2^{-k}+{\LL\choose 2}2^{-\frac{3k}{2}}\right)^{n-\LL(1+\frac{k}{2})-\frac{k}{2}}.
 $$
We bound
 $$
 1-\LL\,2^{-k}+{\LL\choose 2}\,2^{-\frac{3k}{2}}\leq
(1-2^{-k})^\LL\, (1+n^{-1.1})
 $$
 and
 $$
 n-\LL\left(1+\frac{k}{2}\right)-\frac{k}{2} \leq n-k-1
 $$
so that
 $$
 \Pr[\,\OK\,] \leq p^\LL(1+n^{-1.1})^n \leq p^\LL(1+o(1))
 $$
 Our saving comes from
 $$\Pr[z_1\equiv_S z_2]=\Pr[z_3\equiv_S z_4] = 2^{-\LL}.
 $$
The adjacencies on the $z_i$ to $S$ are independent of the event $\OK$.
But
 $$\wedge_{i=1}^\LL I(U_i,y_i)=1\ \ \Rightarrow\  \ \OK.$$
Thus
 $$
 \Pr\left[ \left(\wedge_{i=1}^\LL I(U_i,y_i)\right) \wedge (z_1\equiv_S z_2) \wedge
(z_3\equiv_S z_4)\right] \leq p^\LL\,2^{-2\LL}(1+o(1)).$$

 Putting this together
 $$
 E[Y]  \leq \frac{M^\LL}{\LL!}n^4p^\LL2^{-2\LL}(1+o(1)).
 $$
 Recall that $Mp=\mu \approx 10\log_2 n$.  The function $\mu^x/x!$ hits a maximum
at $x=\mu$ where it is less than $\me^{\mu}$.  Thus
 $$
 \frac{(Mp)^\LL}{\LL!} \leq \me^{\mu}.
 $$
 Hence
 $$
 E[Y] \leq \me^{\mu}n^42^{-2\LL}.
 $$
 We have selected $\LL\approx \mu$ so that 
 $$
  \me^{\mu}2^{-2\LL}= (\me/4)^{\mu(1+o(1))} = n^{-K(1+o(1))},
 $$
 where $K= -10\log_2(\me/4) > 4$.
 \comment{Note that had we chosen $\mu\approx 1000\log_2 n$ we would
have had plenty of room here.  We did need that $\frac{\me}{4}<1$ which
in turn was why we required that {\em two} pair of $z$s be
equivalent.)}%
 We deduce
 $$ E[Y] = o(1)$$
so that almost surely there is no such $\LL$-tuple.
 \comment{For \rit{t4} the
argument is identical except that there are only $n^3$ choices of
$z_1,z_2,z_3$.}%
 Recall that $X$ was the total number of $(U,y)$ with $I(U,y)=1$.
As $E[X]= \mu$
and, from (\ref{var}),
$\Var[X]= O(\mu)$ with probability $1-o(1)$ we have $X\geq \LL$.  
Further, Lemmas~\rlm{U1U2} and~\rlm{6} give that whp the
extensions have properties \rit{t1}--\rit{t2}.  Thus whp
there exists a family of $(U_i,y_i)$ of size $\LL$ which
satisfies \rit{t1}--\rit{t2}.  But also whp
any such family of size $\LL$ will satisfy \rit{t3} and \rit{t4}.
So whp there is such a family.  The expansion of
the family to all $(U,y)$ with $I(U,y)=1$ retains the properties
\rit{t3}--\rit{t4} as the set $S$ is just getting larger.  So
the theorem is proved.\epf

\subsection{Putting All Together}\label{together}

We can now finish the proof of Theorem~\rth{1/2}. By
Lemma~\rlm{AlmostThere} we have whp that all $A\cup
\CY_{k/2}(A)$-similarity classes are singletons except possibly one
$2$-element class $\{x,y\}$. If we let $W=A\cup\{x\}$ and $u=\frac
k2$, then clearly $G$ satisfies all the assumptions of
Lemma~\rlm{1/2}, which implies that $D_2(G)\le k+5$, giving the
required by~\req{k}.

Finally, let us justify the Remark after Theorem~\rth{1/2}. Recall
that given $k$ we have chosen $n$ so that $f(n,k)\approx 10 \log_2 n$
and deduced that $D_2(G)\le k+5$ whp. The probability that the
$(k-1)$-extension property fails for $G$ is at most
 $$
 \binom n{k-1}(1-2^{-k+1})^{n-k}=\me^{k\ln n-2^{-k+1}n-(1+o(1))\, k\ln
k} = f(n,k)\, \me^{-(\frac12+o(1)) k\ln k} = o(1).
 $$
 By Lemma~\rlm{alice}, $k+1\le D(G)$. Thus, $D(G)$ and $D_2(G)$ are
concentrated on at most $5$ different values.

\section{Sparse Random Graphs}\label{sparse}

The following lemma helps us to deal with very sparse random
graphs. Let $t_k=t_k(G)$ be the number of components of $G$ which are
order-$k$ trees. (Thus $t_1(G)$ is the number of isolated vertices.)
For a graph $F$, let $c_F(G)$ be the number of components isomorphic
to $F$.

\blm{comps} Suppose that for any connectivity component $F$ of a graph
$G$ we have
 \beq{comps}
 c_F(G)+v(F)\le t_1(G)+1.
 \eeq
 Then $D(G)=D_1(G)=t_1(G)+2$ unless $G$ is an empty graph (when
$D(G)=D_0(G)=v(G)+1$).
 \elm
 \bpf Assume that $e(G)\not=0$.

The lower bound on $D(G)$ follows by considering $G'$ which is
obtained from $G$ by adding an isolated vertex. The graphs $G$ and
$G'$ are isomorphic as far as non-isolated vertices are concerned. The
best strategy for Spoiler is to pick $t_1(G)+1$ isolated vertices in
$G'$ and, by making one more move in $G$, to show that at least one of
the Duplicator's responses is not an isolated vertex.
  \comment{If $G$ is empty, the last move is unnecessary.}%

On the other hand, let $G'\not\cong G$. There must be a connected
graph $F$ such that $c_F(G)\not= c_F(G')$, say $c_F(G)<
c_F(G')$. Spoiler picks one vertex from some $c_F(G)+1$ $F$-components
of $G'$. If a move $x$ of Duplicator falls into the same component of
$G$ as some her previous move $y$, then Spoiler switches to $G$ and
begins claiming a contiguous path from $x$ to $y$; he wins in at most
$v(F)$ moves by either connecting $x$ to $y$ or by claiming a path of
length $v(F)+1$.

Otherwise, Duplicator must have selected a vertex inside a component
$C$ of $G$ which is not isomorphic to $F$. As soon as this happens,
Spoiler wins by growing a connected set inside the larger component of
the two, in at most $v(F)$ moves.

The total number of moves does not exceed $(c_F(G)+1) + v(F)\le
t_1(G)+2$ (while we have only one alternation), as required.\epf

\bth{1.19} Let $\e>0$ be fixed. Let $p=p(n)\le (\alpha-\e)
n^{-1}$, where $\alpha=1.1918...$ is the (unique) positive root of the
equation
 $$
 \me^{-\alpha +\alpha\me^{-\alpha}} +1 =
 \me^{\alpha\me^{-\alpha}}
 $$
 Then whp $G\in\CG(n,p)$ satisfies the condition~\req{comps}.

 In this range, whp $D(G)=D_1(G)=(\me^{-pn}+o(1))\, n$.
 \eth
 \bpf
 It is easy to compute the expectation of $t_k(G)$ for $G\in\CG(n,p)$:
 $$
 \lambda_k = E[t_k] = {n\choose k}k^{k-2}p^{k-1} q^{k(n-k)+{k\choose
2}-k+1}.
 $$
 Let $c=pn$. For a fixed $k$ we have $\lambda_k=(f_k+o(1)) n$, where
$f_k:=\frac{c^{k-1} k^{k-2}}{k!\, \me^{ck}}$. We have
 $$
 \frac{f_{k+1}}{f_k}
  = c\me^{-c} \times (1+1/k)^{k-2}.
 $$
 The first factor $c\me^{-c}$ is at most $1/\me$ (maximized for
$c=1$). Unexciting algebraic calculations show that the second factor
is monotone increasing for $k\ge 1$ and approaches $\me$ in the
limit. This implies that the sequence $f_k$ is decreasing in $k$. (In
particular, $f_1$ is strictly bigger than any other $f_i$, $i\ge 2$.)

Also, $f_1=\me^{-c} > 0.3$ for $0\le c\le \alpha$. 

Theorem~5.7 in Bollob\'as~\cite{bollobas:rg} describes the structure
of a typical $G$ for $p=O(n^{-1})$. In particular, it implies that there
is a constant $K$ such that whp at least $0.9 n$ vertices of $G$
belong either to tree components of orders at most $K$ or to the giant
component.
 \comment{Bela does not refer to the origins of the result. Whose
result is it?}%
 The giant component (for $c>1$) has order $(1-\frac sc+o(1))\, n$,
where $s$ is the only solution of $s\me^{-s}=c\me^{-c}$ in the range
$0<s\le 1$. It is routine to see that $f_1 > 1-\frac sc$. (In fact,
$c=\alpha$ is the root of $f_1=\frac sc$; this is where $\alpha$ comes
from.)

A theorem of Barbour~\cite{barbour:82} (Theorem~5.6
in~\cite{bollobas:rg}) implies that, for any $k\le K$, we have whp
 \beq{barbour}
 |t_k(G) - \lambda_k|\le o(n).
 \eeq

Now, we have all the ingredients we need to check~\req{comps}. Let
$F\subset G$ be any connectivity component. If $F$ is the giant
component of $G$, then $c_F(G)=1$ but, as we have seen, $v(F)\le
t_1(G)$ so~\req{comps} holds. So we can assume that $v(F)=o(n)$. If
$v(F)> K$, then $c_F(G)+v(F)\le \frac{0.1 n}K +K < t_1(G)$. If $F$ is a tree
with $k\in[2,K]$ vertices, then~\req{barbour} and the inequality $f_1>
f_k$ imply the required. Finally, it remains to assume that the
component $F$ of order at most $K$ contains a cycle. But the expected
number of such components is at most $\sum_{i=2}^K {n\choose i} p^{i}
i^{i-1} {i\choose 2} = O(1)$. Markov's inequality implies that whp no
such $F$ violates~\req{comps}.\epf

Of course, the value of $t_1(G)$ can be estimated more precisely for
some $p$ than we did in Theorem~\rth{1.19}.  Without going into much
details, let us describe some of the cases here. Let $\omega$ be any
function of $n$ which (arbitrarily slowly) tends to the infinity with
$n$.

If $n^2p\to 0$, then whp we have isolated vertices and edges only. The
distribution of $t_2(G)=e(G)$ approaches the Poisson distribution
$\Poisson{\lambda_2}$. Hence, we have whp that $n-\omega< D(G)\le
n+1$.

Suppose that $n^2p\not\to 0$ but $pn\to 0$. The expected number of
vertices in components of order at least $3$ is at most $n {n\choose2}
3p^2=o(\lambda_2)$. By Markov's inequality, whp we have $o(\lambda_2)$
such vertices. On the other hand, the distribution of $t_2(G)$ is
$o(1)$-close to $\Poisson{\lambda_2}$ (Theorem~5.1
in~\cite{bollobas:rg}). Hence,
 $$
 \Big| D(G)- n+ \frac{\lambda_2}2\Big| \le o(\lambda_2)+\omega.
 $$

Observe that there is no phase transition in the behavior of $D(G)$ at
$p\approx \frac1n$. This should not be surprising: $D(G)$ is
determined by $t_1(G)$ in this range. We believe (but were not able to
prove) that whp $D(G)= t_1(G)+2$ for $p=O(n^{-1})$. To show this it is
enough to define the giant component by a sentence of depth $o(n)$. In
fact, we conjecture that a far stronger claim is true.

\bcj{giant} Let $\frac{1+\e}n< p=O(n^{-1})$. Then whp $D(F)=O(\ln n)$,
where $F$ is the giant component of $G\in\CG(n,p)$.
 \ecj

\section{Modeling Arithmetics on Graphs}

In this section we consider $D(G)$ for the random graph $G\in
\CG(n,p)$ where $p=n^{-1/4}$.  We expect that our results would hold
for $p=n^{-\ah}$ for any {\em rational} $\ah\in (0,1)$, but this would
require considerable technical work so we are content with this one
case.  In~\cite[Section 8]{spencer:slrg} it was shown, for
$\ah=\frac{1}{3}$, that there was an arithmetization of certain sets
that led to non-convergence and non-separability results.  Our methods
here will be similar.

\bth{logstar} Let $p=n^{-1/4}$ and $G\in\CG(n,p)$.  Then whp
\[  D(G) = \Theta(\log^*n).  \]
 \eth

The lower bound is very general.  We use only the trivial fact that
any particular graph is the value of the random graph with probability
$O(n^{-1})$.  (Indeed, the value is exponentially small.) Let $F(k)$
be the number of first order sentences of depth at most
$k$.  Then $\Pr[D(G)\leq k]= O(F(k)n^{-1})$ as there are at most
$F(k)$ such graphs.  From general principles~\cite[Section
2.2]{spencer:slrg} $F(k)$ is bounded effectively by the tower function
so that if $k=c\log^*n$ with $c$ appropriately small $F(k)=o(n)$ and
$\Pr[D(G)\leq k]=o(1)$.

Now we turn to the main part, bounding $D(G)$ from above.  For any set
$W$ of vertices let $N(W)$ denote the set of common neighbors of $W$.
When $|W|=4$, $\Pr[N(W)=\emptyset] = (1-p^4)^{n-4} \approx \me^{-1}$.  We
are guided by the idea that $N(W)=\emptyset$ is like a random symmetric $4$-ary
predicate with probability $\me^{-1}$, which is bounded away from both zero and
one.

Let $W$ be a set of four vertices.  Dependent only on $W$ we define
 \begin{itemize}
 \item  $A=N(W)$, the common neighbors of $W$;
 \item $B$, consisting of those $z\not\in W\cup A$ such that $z$ is adjacent
to precisely four vertices of $A$ and no other $z'\not\in W\cup A$ has
exactly the same adjacencies to $A$.
 \end{itemize}

For $w\not\in W\cup A$ let $H_w(A)$ denote the $3$-regular hypergraph
on $A$ consisting of those triples $T$ so that there is no $z\not\in
W\cup A$ adjacent to $T\cup \{w\}$.  (The condition that $z\not\in
W\cup A$ is a technical convenience that does not asymptotically
affect the $H_w(A)$.)  If, further, $a\in A$ we let $H_{w,a}(A)$
denote the $2$-regular hypergraph (i.e.\ graph) of pairs $T$ with
$T\cup\{a\}\in H_w(A)$.  Further, for distinct $a,b\in A$ we let
$H_{w,a,b}(A)$ denote the $1$-regular hypergraph (i.e.\ set) of
elements $y$ with $\{a,b,y\}\in H_w(A)$.  For $w\not\in W\cup A\cup B$
let $H_w(B)$ denote the $3$-regular hypergraph on $B$ consisting of
those triples $T$ so that there is no $z\not\in W\cup A\cup B$
adjacent to $T\cup \{w\}$.  (Again, the condition $z\not\in W\cup
A\cup B$ is a technical convenience.)  Informally, the idea is that
the $H_w(A), H_w(B)$ act like random objects with probability
$\me^{-1}$.

Call $A$ universal if the $H_v(A)$, $v\not\in W\cup A$,  range over all $3$-regular 
hypergraphs on $A$.  Call $B$ splitting if the $H_v(B)$, $v\not\in W\cup A\cup B$, 
are all different.  
As there are $2^{\Theta(m^3)}$ $3$-regular hypergraphs over an $m$-set a simple
counting argument gives that if $A$ is universal we must have $|A|=O(\ln^{1/3}n)$
while if $B$ is splitting we must have $|B|=\Omega(\ln^{1/3}n)$.

Our argument splits into two lemmas.  

\blm{1}  Whp there exists a $4$-set $W$  such that, with
$A,B$ as defined above,
 \ben
\item $A$ is universal; 
\item $B$ is splitting.
 \een
\elm

\blm{2}  Any graph $G$ on $n$ vertices with the property of Lemma~\rlm{1} has
$D(G)= O(\log^*n)$.\elm

Note that the proof of Lemma~\rlm{1} is a random graph argument while
the proof of Lemma~\rlm{2} is a logic argument involving no
probability.

\smallskip\noindent\textit{Proof of Lemma~\rlm{1}.}  Set
$u=\lfloor\ln^{0.3}n\rfloor$.  For any set $W$ of four vertices
 $$
 \Pr[\, |N(W)|=u\,]= \Pr[\,\mathrm{Bin}(n-4,p^4)=u\,] \approx \me^{-1}/u!.
 $$
  Thus the expected number $\mu$ of such $W$ has $\mu\approx {n\choose
4}\me^{-1}/u!$ which approaches infinity. An elementary second moment
calculation gives that the number of such $W$ is $(1+o(1))\,\mu$ whp.
Hence it suffices to show that the expected number of $W$ with $A$
having size $u$ but $A,B$ failing the conditions of Lemma~\rlm{1} is
$o(\mu)$.  Fix $W$ and $A$ of size $u$.  It suffices to show that
$A,B$ satisfy Lemma~\rlm{1} whp.  The conditioning is only on the
adjacencies involving a vertex of $W$, all other adjacencies remain
random.

First we show that $A$ is universal.  Let $Z$ be those vertices
adjacent to four or more vertices of $A$.  Let $Z'$ consist of
vertices with at least one neighbor in $Z$. Whp every four vertices in
the graph have $O(\ln n)$ common neighbors so $|Z|$ is polylog
while $|Z'|=n^{3/4+o(1)}$.

For each $3$-set $Y\subset A$ let $N^-(Y)$ denote those $v\not\in
W\cup A\cup Z$ which are adjacent to all vertices of $Y$.  The $|N^-(Y)|$
are independent, each with Binomial distribution
$\mathrm{Bin}(n-o(n),(1+o(1))\,n^{-3/4})$.  The probability that
$|N^-(Y)|>2n^{1/4}$ or $|N^-(Y)|<\frac{1}{2}n^{1/4}$ is then less than
$\exp(-cn^{1/4})$ for a constant $c$ by Chernoff bounds.  We only need
that this probability is $o(n^{-3})$.  Thus with high probability for
every $3$-set $Y\subset A$ we have $|N^-(Y)|\in
[\frac{1}{2}n^{1/4},2n^{1/4}]$.  Condition on these $N^-(Y)$ satisfying
these conditions.  Let $R$ be the (remaining) vertices, not in $W,A,Z,Z'$
nor any of the $N^-(Y)$.  For $z\in R$ the adjacencies to the $N^-(Y)$ are
still random.  For such $z$ we have $Y\in H_z(A)$ if and only if $z$
is adjacent to no vertex in $N^-(Y)$. (Note that $z$ does not send any
edges to $Z$.) As $|N^-(Y)|\leq 2n^{1/4}$ the
probability that $z$ is adjacent to no vertex of $N^-(Y)$ is at least
$\me^{-2}$.  As $|N^-(Y)|\geq \frac{1}{2}n^{1/4}$ the probability that
$z$ is adjacent to some vertex of $N^-(Y)$ is at least $1-\me^{-1/2}$.
Set $\gam=\min(\me^{-2},1-\me^{-1/2})$.  Then for any hypergraph $H$
on $A$ we have $\Pr[H_z(A)=H]\geq \gam^{{u\choose 3}}$.  But these
events are now independent over the $(1-o(1))\,n$ values $z\in R$ so
that the probability that no $H_z(A)=H$ is less than
$(1-\gam^{{u\choose 3}})^n$.
Here because $u^3=o(\ln n)$ this quantity is less than,
say, $\exp(-n^{0.99})$.  There are fewer than $2^{u^3}$ hypergraphs
$H$ on A.  Hence the probability that any such $H$ is not one of the
$H_z(A)$ is less than $2^{u^3}\exp(-n^{0.99})$.  The $2^{u^3}$ term is
basically negligible and the probability that $A$ is not universal is
less than $\exp(-n^{0.98})$ and certainly $o(1)$.  We note that $A$
being universal will not be fully needed in Lemma~\rlm{2}, we shall need
only seven particular values of $H_z(A)$.

Now we look at the size of $B$. For each $z\not\in A\cup W$ the
probability that $z$ is adjacent to precisely four elements of $A$ is
${u\choose 4}p^4(1-p)^{u-4} \approx u^4n^{-1}/24$ and given this the
probability that no other $z'$ has the same adjacencies is approximately
$(1-p^4)^n\approx \me^{-1}$ so $B$ has expected size $\mu
\approx u^4/24\me$.  A second moment calculation gives that whp
$|B|\approx \mu= \Theta(\ln^{1.2}n)$.

Finally we show that $B$ is splitting.  At this stage $W,A,B$ are fixed and
all of the adjacencies that do not have at least one vertex from $W\cup A$ are
random.  Whp no $z\not\in W\cup A\cup B$ is adjacent to five (or more) vertices
of $B$.  Let $Z$ be those $z\not\in W\cup A\cup B$ adjacent to four vertices
of $B$.  Whp $|Z|$ is polylog.  

For each $3$-set $Y\subset B$ we let $N^*(Y)$ 
denote those $v\not\in W\cup A\cup B$ which are adjacent to all vertices
of $Y$  and $N^-(Y)=N^*(Y)-Z$.
As before, whp all $|N^*(Y)|$ have size between $\frac{1}{2}n^{1/4}$ and
$2n^{1/4}$ and so the same, asymptotically, holds for the $|N^-(Y)|$.
As before we fix the $N^*(Y)$ and their adjacencies to $Y$. Consider distinct 
$u,u'\not\in W\cup A\cup B$.  The probability that either $u$ or $u'$
is adjacent to nine (or more) vertices of $Z$ is $o(n^{-2})$.  Call a
$3$-set $Y\subset B$ exceptional if $u$ or $u'$ is adjacent to some $z\in Z$ 
which is adjacent to all of $Y$.  With probability $1-o(n^{-2})$ there
are at most $2\cdot 8\cdot 4 = 64$ exceptional $Y$.  Hence the number
of nonexceptional $Y$ is $\approx {|B|\choose 3}$.
For the nonexceptional
$Y$ we have $Y\in H_{u}(B)$ if and only if $u$ is adjacent to no vertex
in $N^-(Y)$ and similarly for $Y\in H_{u'}(B)$.  Thus $\Pr[Y\in H_{u}(B)]\in
[\gam,1-\gam]$ with $\gam$ as previously defined.  Further these events are
independent over different nonexceptional $Y$.  Set
$\gam^*=\gam^2+(1-\gam)^2$.  Then $H_u(B), H_{u'}(B)$ agree on a
nonexceptional $Y$ with probability at most $\gam^*$.  Independence
gives that they agree on all nonexceptional $Y$ with probability at
most $\gam^*$ to the  $\approx {|B|\choose 3}$ power. As this power
is $\gg \ln n$ the probability is certainly $o(n^{-2})$.  There are $O(n^2)$
choices of $u,u'$ so whp no  $H_u(B)=H_{u'}(B)$.\epf

\smallskip\noindent\textit{Proof of Lemma~\rlm{2}.}  The main portion of the argument consists of placing an
arithmetic structure on $A$ in such a way that any vertex in $A$ can be
described with quantifier depth $O(\log^*|A|)=O(\log^*n)$.  For convenience
we assume $|A|=3s+2$.  (Otherwise we would give the one or two extra elements
of $A$ which would increase the depth by at most two.)
Label the elements of $A$
by $a,b,x_1,\ldots,x_s,y_1,\ldots,y_s,z_1,\ldots,z_s$
in an arbitrary way.  Now we, effectively,  model arithmetic on $A$.  
From the universality there exist $w_1,\ldots,w_7$ (witnesses) such that
 \ben
\item $H_{w_1}$ consists of the triples $\{x_i,y_i,z_i\}$
\item $H_{w_2,a,b}$ consists of the elements $x_i$
\item $H_{w_3,a,b}$ consists of the elements $y_i$
\item $H_{w_4,a}$ consists of all pairs $\{x_i,y_j\}$ with $i\leq j$
\item $H_{w_5}$ consists of all triples $\{x_i,y_j,z_{i+j}\}$ with $1\leq i,j,i+j \leq s$ 
\item $H_{w_6}$ consists of all triples $\{x_i,y_j,z_{i\cdot j}\}$ 
with $1\leq i,j,i\cdot j \leq s$.
\item $H_{w_7,a}$ consists of all pairs $\{x_i,y_{2^i}\}$ with $1\leq i\leq s$ and $2^i\leq
s$.
\een
We give a first order expression in terms of the $v_1,v_2,v_3,v_4$ which define $A$
and the special elements $a,b\in A$ and the witnesses $w_1,\ldots,w_7$ which forces
$A$ to have this form.  Note that membership in $A$ is given by a first order statement
and membership in an $H_w$ or $H_{w,a}$ or $H_{w,a,b}$ is given by a first order
statement in terms of the variables.  Let $A^-=A-\{a,b\}$ for convenience.
Now we express the following six properties.
 \ben
\item ($1$-Factor) $H_{w_1}$ consists of vertex disjoint triples and every element
of $A^-$, and only those elements, are in such a triple.
\item (Splitting the $1$-Factor) For each triple in $H_{w_1}$ exactly one of the
elements is in $H_{w_2,a,b}$ and exactly one (a different one) is in $H_{w_3,a,b}$.
Now let
(for convenience)
$X$ denote those $x\in A^-$ with $x\in H_{w_2,a,b}$ and let $Y$ denote those
$y\in A^-$ with $y\in H_{w_3,a,b}$ and $Z$ the other elements of $A^-$.  Henceforth
the use of the letter $x,y,z$ shall tacitly assume that the element is in the
respective set $X,Y,Z$.  We write $x\sim y$ or $x\sim z$ or $y\sim z$ if the two
elements are in a common triple in $H_{w_1}$.
\item (Creating $\leq$) Here adjacency is in $H_{w_4,a}$.  We require that all
adjacencies be between an $x$ and a $y$.
Let $N(x)$ denote the
$y$ adjacent to $x$.  We require that for every $x,x'$ either $N(x)\subset N(x')$
or $N(x')\subset N(x)$ with equality only when $x=x'$.  We require that when $y\sim x$
then $y\in N(x)$.  This forces the $N(x)$ to form a chain and so the $x$ and $y$ can
be renumbered to fit the condition.  We now define $x\leq x'$ by $N(x')\subseteq N(x)$.
The relations $\geq,>,<$ have their natural Boolean meaning in terms of $\leq$.
We define $y\leq y'$ and $z\leq z'$ by $x\leq x'$ where $x\sim y\sim z$ and $x'\sim y'
\sim z'$.  We let $x_1,y_1,z_1$ denote the first elements under $\leq$ and
$x_s,y_s,z_s$ the last elements.  The notions of successor $x^+$ and predecessor
$x^-$ are naturally defined (when they exist) in terms of $\leq$. We let $z_2$ denote
the successor of $z_1$.
\item (Creating addition) Addition is generated from the formulas $\ah+1=\ah^+$ and
$\ah+\beta^+=(\ah+\beta)^+$, though we need some care as addition in this model is
not always defined.  
For every $x\in X, y\in Y$ there is at most one $z\in Z$
with $\{x,y,z\}\in H_{w_5}$. $\{x_1,y_1,z_2\}\in H_{w_5}$.  For
$x\neq x_s$, $\{x,y_1,z\}\in H_{w_5}$ where $z\sim x^+$.  
If $\{x,y,z\}\in H_{w_5}$ 
and $y,z$ have successors then $\{x,y^+,z^+\}\in H_{w_5}$.  
If $\{x,y,z\}\in H_{w_5}$ and $y,z$ have predecessors then
$\{x,y^-,z^-\}\in  H_{w_5}$.  We let $x+x'=x^*$ denote
that $\{x,y',z^*\}\in H_{w_5}$ where $y'\sim x'$ and $y^*\sim x^*$.
Let $x+z=z'$ mean that when $z,z'$ are replaced by their $\sim$ elements in
$x$ that then we have the equality, and similarly for other forms like $y+y'=z$.
\item (Creating multiplication) Multiplication is generated from the formulas 
$\ah\cdot 1=\ah$ and
$\ah\cdot\beta^+=(\ah\cdot\beta)+\beta$, though we need some care as addition in this model is
not always defined.
For every $x\in X, y\in Y$ there is at most one $z\in Z$
with $\{x,y,z\}\in H_{w_6}$. $\{x,y_1,z\}\in H_{w_5}$ precisely when $x\sim z$.
If $\{x,y,z\}\in H_{w_6}$ and $y$ has a successor then $\{x,y^+,z'\}\in H_{w_6}$
if and only if $z'=x+z$.
\item (Creating exponentiation) Base two exponentiation is defined by $2^1=2$ and
$2^{\ah^+}=2^{\ah}+2^{\ah}$, though we need some care as addition in this model is
not always defined.
For every $x\in X$ there is at most one  $y\in Y$
with $\{x,y\}\in H_{w_7,a}$.  $\{x_1,y_2\}\in H_{w_7,a}$.  If $\{x,y\}\in H_{w_7,a}$ 
then $\{x^+,y'\}\in H_{w_7,a}$ if and only if $y+y=y'$.  We write $x'=2^x$ if
$\{x,y'\}\in H_{w_7,a}$ and $x'\sim y'$.
\een

We can give the $d$-th binary digit of $x$ (we count the first digit
as the first on the right, zero if and only if $x$ is even) for all
$x\in X$. The $1$-st digit of $x$ is zero if and only if $x=x'+x'$ for
some $x'$.  Otherwise the $d$-th digit is zero if and only if there
exist $q,r$ $x=q\cdot 2^d+r$ with $r<2^{d-1}$ or if $x=q\cdot 2^d$ or
if $x<2^{d-1}$.  (This technical complication is caused by leaving
zero out of the model.)  This is first order as we already have
multiplication, exponentiation, less than and addition.

Now any $x\leq s$ is described with quantifier depth $\Theta(\log^*s)$.
Let $x$ have $m$ digits. We say that $x\leq 2^m$ and
the disjunction over $d\leq m$ of the statements that the $d$-th binary digit
is what it is.  For each such $d$ (and for $m$) we have to describe $d$.  But now we are
describing numbers up to $\log_2 s$
and so by an induction the total
depth will be $\Theta(\log^*s)$.  This also includes describing the last element
$x_s$ so that we determine $s$ with depth $\Theta(\log^*s)$.

Now the elements $v$ of $B$ are described with depth
$\Theta(\log^*s)$ by describing the four vertices of $A$ that $v$ is
adjacent to and saying that $v$ is adjacent to no other vertices of $A$
and that no other $w\not\in W\cup A$ has just those adjacencies to
$A$.

Finally any $v\not\in W\cup A\cup B$ is described in depth $O(\log^*
s)$ by listing the edges of $H_v(B)$ and stating that no other $v'$
produces the same hypergraph.  The assumption $B$ splitting means that
we have described all the vertices.\epf

\section{Concluding Remarks}\label{open}

Our Theorem~\rth{constant} has a strong link with the zero-one law
which was discovered independently by Glebskii et al \cite{glebskii+:69} and
Fagin \cite{fagin:76} and says that $G\in\CG(n,\frac12)$
satisfies any fixed first order sentence with probability approaching
either 0 or~1.  Given $\epsilon\in(0,1)$, define $T_\epsilon(n)$ to be
the maximum $k$ such that, if $n_1,n_2\ge n$, then $D(G,H) > k$ with
probability at least $1-\epsilon$ for independent $G\in\CG(n_1,p)$ and
$H\in\CG(n_2,p)$. The Bridge Theorem \cite[Theorem
2.5.1]{spencer:slrg} says that, in a rather general setting, a
zero-one law is obeyed iff, for each $\epsilon$, $T_\epsilon(n)$ tends
to the infinity as $n$ increases.  Spencer and St.~John
\cite{spencer+stjohn:01} call $T_\epsilon(n)$ the \emph{tenacity
function} and suggest it as a quantitative measure for observation of
a zero-one law.  While in \cite{spencer+stjohn:01} the tenacity
function is studied for words, here we are able to find its
asymptotics in the case of graphs. Since the lower bound based on the
$k$-extension property goes through for $D(G,H)$ with both $G$ and $H$
random, we have $T_\epsilon(n)=\log_{1/p} n + O(\ln \ln n)$,
irrespective of the constant $\epsilon$.

Another interesting first order parameter of a graph $G$ is $I(G)$,
the smallest depth of a sentence distinguishing $G$ from any
non-isomorphic graph of the \emph{same} order as $G$. Of course,
$I(G)\le D(G)$, so all upper bounds we have proved apply to $I(G)$ as
well. All our lower bounds also apply to $I(G)$ with the exceptions of
Theorem~\rth{1.19}. Its $I(G)$-analog would say that, for
$G\in\CG(n,\frac{c}n)$ with $c< c_0$, where $c_0=1.034...$, we have
whp
 \beq{SparseI}
 I(G)= (1+o(1))\, t_2(G) = (c\me^{-2c}/2+o(1))\, n,
 \eeq
 where $t_2(G)$ denotes the number of isolated edges. The reason is
that if $G\not\cong G'$ but $v(G)=v(G')$, then the multiplicities of
at least two non-isomorphic components must differ while the two most
frequent components in $G$ are isolated vertices and edges. (And the
order of the giant component catches up with $t_2(G)$ at $p\approx
\frac{c_0}n$.)  The Reader should not have any problem in filling up
the missing details.

We make the following general conjecture.

\bcj{general} Let $\e>0$ be fixed and $n^{-1+\e}\le p \le \frac12$. Then
whp $D(G)=O(\ln n)$.\ecj

One can also ask about $D_\#$, the analog of $D(G)$ when we add
counting to first order logic. Here, the situation is strikingly
different. A result of Babai and Ku\v cera~\cite{babai+kucera:79}
(combined with Immerman and Lander's~\cite{immerman+lander:90} logical
characterization of the vertex refinement step
in~\cite{babai+kucera:79}) implies that whp $G\in\CG(n,\frac12)$ can
be defined by a first order sentence with counting of quantifier depth
at most~$4$.

\OPend{misc,graph,random,general,ramsey}